\documentclass[12pt]{article}
\usepackage{amsmath, amssymb, amsthm, amsfonts}

\def\stir#1#2{\genfrac{\{}{\}}{0pt}{}{#1}{#2}}
\def\hot{+ \text{higher-order terms}}
\def\RR{\mathbb{R}}

\newtheorem{theorem}{Theorem}

\title{Another combinatorial determinant}
\author{Kiran S. Kedlaya \\ Massachusetts Institute of Technology}
\date{May 18, 1999}

\begin{document}
\maketitle

The following result, a special case of a theorem of Mina \cite{mina},
was recently given an elegant proof by
Wilf \cite{wilf}.
\begin{theorem}
Let $f = 1 + a_1 x + a_2 x^2 + \cdots$ be a formal power series, and
define a matrix $c$ by
\[
c_{i,j} = [x^j] f^i \qquad (i,j \geq 0).
\]
(Here $[x^j] f$ denotes the coefficient of $x^j$ in $f$.)
Then
\[
\det \left( (c_{i,j})_{i,j=0}^n\right)  = a_1^{n(n+1)/2} \qquad
(n = 0,1,2,\dots).
\]
\end{theorem}
The purpose of this paper is to show that this result remains essentially
unchanged if we take powers in the sense of composition, instead of
multiplication.
\begin{theorem}
Let $f = x + b_1 x^2 + b_2 x^3 + \cdots$ be a formal power series, and define
$f^{(0)} = x$ and $f^{(i)} = f(f^{(i-1)})$ for $i>0$. Define a matrix $c$ by
\[
c_{i,j} = [x^{j+1}] f^{(i)} \qquad (i,j \geq 0).
\]
Then
\[
\det \left( (c_{i,j})_{i,j=0}^n\right)  = 1! 2! \cdots n! b_1^{n(n+1)/2} \qquad
(n = 0,1,2,\dots).
\]
\end{theorem}
In fact, we prove both of these theorems at once, by formulating and
proving a common generalization. In both theorems, each row of the matrix is
obtained from the previous row
by applying a certain transformation of power series: in
Theorem~1, the transformation is
\[
t \mapsto t(1 + a_1 x + a_2 x^2 + \cdots),
\]
while in Theorem~2, the transformation is
\[
t \mapsto t + b_1 t^2 + b_2 t^3 + \cdots.
\]
This suggests that more generally, we should consider transformations
of the form
\[
t \mapsto f(t) = \sum_{m=1}^\infty \sum_{n=0}^\infty b_{m,n} t^{m} x^n.
\]
\begin{theorem}
Let $f(t) = \sum_{m=1}^\infty \sum_{n=0}^\infty b_{m,n} t^{m} x^n$
be a formal power series in two variables $t$ and $x$, and assume
$b_{1,0} = 1$.
Define $f^{(0)}(t) = t$ and $f^{(i)}(t)  = f(f^{(i-1)}(t))$ for $i > 0$.
Define a matrix $c$ by
\[
c_{i,j} = [x^{j+1}] f^{(i)}(x) \qquad (i,j \geq 0).
\]
Then
\[
\det \left( (c_{i,j})_{i,j=0}^n\right)
= \prod_{k=1}^n
\left( \sum_{m=0}^k m! \stir{k+1}{m+1} b_{2,0}^m b_{1,1}^{k-m} \right)
\qquad (n = 0,1,2,\dots),
\]
where $\stir{x}{y}$ is the Stirling number of the second kind (the
number of partitions of $x$ labeled objects into $y$ nonempty sets).
\end{theorem}
\begin{proof}
Following \cite{wilf}, we prove that the matrix
\[
b_{i,j} = (-1)^{i+j} \binom{i}{j} \qquad (i,j \geq 0)
\]
has the property that $bc$ is upper triangular. Put
\[
g_i(t) = \sum_k (-1)^{i+k} \binom{i}{k} f^{(k)}(t);
\]
then $(bc)_{i,j} = [x^j] g_i(x)$, and the theorem will follow from
the fact that
\[
g_i(t) = \sum_{j=0}^i j! \stir{i+1}{j+1} b_{2,0}^{j} b_{1,1}^{i-j} t^{j+1} x^{i-j} \hot,
\]
which we prove by induction on $i$ (the case $i=0$ being true by
definition).
If we write $g_{i-1}(t) = \sum_{m,n} k_{m,n} t^{m+1} x^n$ (so in particular,
$k_{m,n} = 0$ for $m+n < i-1$ and $k_{m,n} =
m! \stir{m+n+1}{m+1} b_{2,0}^{m} b_{1,1}^n$ for $m+n = i-1$), then
\begin{align*}
g_i(t) &= g_{i-1}(f(t)) - g_{i-1}(f(t)) \\
&= \sum_{m,n} k_{m,n} x^n (f(t)^{m+1} - t^{m+1}) \\
&= \sum_{m,n} k_{m,n} x^n t^{m+1} [(m+1)(f(t)/t - 1) - \binom{m+1}{2}(f(t)/t - 1)^2 + \cdots] \\
&= \sum_{m+n=i-1} (m+1) k_{m,n} t^{m+1} x^n (b_{2,0} t + b_{1,1} x) \hot \\
&= \sum_{m+n=i} t^{m+1} x^n [(m+1) b_{1,1} k_{m,n-1} + m b_{2,0}
k_{m-1,n}] \hot \\
&= \sum_{m+n=i} t^{m+1} x^n m! b_{2,0}^{m} b_{1,1}^n \left(\stir{i}{m} +
(m+1) \stir{i}{m+1}\right) \hot \\
&= \sum_{m+n=i} t^{m+1} x^n m! b_{2,0}^{m} b_{1,1}^n \stir{i+1}{m} \hot,
\end{align*}
as desired.
\end{proof}

It should be pointed out that the full version of Mina's theorem, which
states that
\[
\det \left\{ \left( \frac{d^j}{dx^j} f(x)^i \right)_{i,j=0}^n
 \right\} = 1!2!\cdots n!f'(x)^{n(n+1)/2},
\]
does not appear to admit an analogous generalization. The difficulty
seems to be that while Mina's theorem follows from applying Theorem~1
to the Taylor expansion of $f(x)/f(t)$ at $t$ for each $t \in \RR$,
Theorem~2 can only
be applied to the Taylor expansion of $f$ at its fixed points.

\end{document}